\documentclass[%
 aip,
 cha,
 sd,%
 amsmath,amssymb,
 preprint,%
]{revtex4-1}

\usepackage{graphicx}% Include figure files
\usepackage{dcolumn}% Align table columns on decimal point
\usepackage{bm}% bold 
\usepackage{algorithm}
\usepackage{subfig}

\begin{document}

\preprint{AIP/}

\title{Extended dynamic mode decomposition with dictionary learning: a data-driven adaptive spectral decomposition of the Koopman operator}
%\thanks{Footnote to title of article.}

\author{Qianxiao Li}\email{liqix@ihpc.a-star.edu.sg.}
\affiliation{Institute of High Performance Computing, Singapore}
\author{Felix Dietrich}%
\affiliation{ 
Faculty of Mathematics, Technical University of Munich, Germany
}%
\author{Erik M. Bollt}%
\affiliation{ 
Department of Mathematics \& Department of Electrical and Computer Engineering, Clarkson University, USA
}%

\author{Ioannis G. Kevrekidis}
\email{yannis@princeton.edu}
\affiliation{%
Department of Chemical and Biological Engineering \& the Program in Applied and Computational Mathematics, Princeton University, USA
}%

\date{\today}% It is always \today, today,

\begin{abstract}
Numerical approximation methods for the Koopman operator have advanced considerably in the last few years. In particular, data-driven approaches such as dynamic mode decomposition (DMD) and its generalization, the extended-DMD (EDMD), are becoming increasingly popular in practical applications.
The EDMD improves upon the classical DMD by the inclusion of a flexible choice of dictionary of observables that spans a finite dimensional subspace on which the Koopman operator can be approximated. This enhances the accuracy of the solution reconstruction and broadens the applicability of the Koopman formalism. Although the convergence of the EDMD has been established, applying the method in practice requires a careful choice of the observables to improve convergence with just a finite number of terms. This is especially difficult for high dimensional and highly nonlinear systems. In this paper, we employ ideas from machine learning to improve upon the EDMD method.
We develop an iterative approximation algorithm which couples the EDMD with a {\it trainable} dictionary represented by an artificial neural network. Using the Duffing oscillator and the Kuramoto Sivashinsky PDE as examples, we show that our algorithm can effectively and efficiently adapt the trainable dictionary to the problem at hand to achieve good reconstruction accuracy without the need to choose a fixed dictionary {\it a priori}. Furthermore, to obtain a given accuracy we require fewer dictionary terms than EDMD with fixed dictionaries. This alleviates an important shortcoming of the EDMD algorithm and enhances the applicability of the Koopman framework to practical problems. 
\end{abstract}

%Use showkeys class option if keyword
\maketitle

\begin{quotation}
Every dynamical system has an associated Koopman operator, which encodes many important properties of the system. Most notably, it characterizes the temporal evolution of observables by linear, albeit infinite-dimensional, dynamics even when the underlying dynamical system is non-linear. In recent years, the growing availability of data and novel numerical techniques have enabled us to study this operator computationally. 
Extended Dynamic Mode Decomposition (EDMD) is one such technique for approximating the spectral properties of the operator. While effective for some problems, a clear drawback of EDMD is the requirement to select \emph{a priori} a suitably efficient collection of basis functions, called a \emph{dictionary}. 
In this paper, we use ideas from machine learning to optimally adapt the dictionary to data. This enables us to obtain improved numerical approximations without resorting to large dictionary sizes. 
We demonstrate the efficiency of our algorithm by approximating the Koopman operator for the Duffing oscillator system and also the Kuramoto-Sivashinsky PDE.
\end{quotation}

\section{Introduction}

In the analysis of dynamical systems, a primary object of study is the state and its evolution. In this traditional setting, powerful tools from differential geometry can characterize dynamical systems by their trajectories and invariant manifolds in phase space. 
In recent years, however, advances in numerical techniques and the broader availability of data have sparked renewed interest in an alternative, {\it operator} view on dynamical systems \cite{mezic-2004,mezic-2005}: the Koopman operator framework \cite{koopman1931hamiltonian}. In this framework, the central objects of study are {\it observables}, which are functions of the state of the dynamical system. The Koopman operator then describes the temporal evolution of these functions driven by the underlying dynamics. 

The Koopman formalism is useful in several ways. First, the Koopman dynamics is linear, albeit infinite-dimensional, and hence amenable to powerful methods from operator theory such as spectral analysis \cite{neumann1932operatorenmethode,halmos1942operator,halmos1957introduction}. Second, it is especially suited for studying high-dimensional systems, where the phase space is so large that little can be said from the differential geometry point of view.
The Koopman approach allows one to focus on the evolution of a lower number of observables. In fact, in applications this is often the case: the evolution of a small number of measurements (observables) of an otherwise high dimensional system is recorded. Lastly, from a numerical point of view, it allows one to employ traditional techniques in numerical linear algebra to perform linearization and ``normal mode analysis'' for nonlinear systems. This is an important advantage of the Koopman framework for current challenges in model reduction, prediction, data fusion, and system control \cite{rowley2009spectral,budisic-2012,williams-2015b,giannakis-2015b,brunton-2016b,korda-2016}. Applications range from fluid dynamics \cite{mezic-2013,sharma-2016}, energy modeling in buildings \cite{georgescu-2015} and oceanography \cite{giannakis-2015b}, to molecular kinetics \cite{wu-2017} and beyond.

We note that, given appropriate function spaces, the adjoint of the Koopman operator is the Perron-Frobenius operator. It operates on phase-space density functions and advances them in time according to the underlying dynamics. The duality of these two operators can be described as ``dynamics of observables'' for the Koopman operator in contrast to ``dynamics of densities'' for the Perron-Frobenius operator \cite{budisic-2012}. Both are valid descriptions of the underlying system through the perspective of linear operators. 

In effect, the Koopman framework converts a finite-dimensional, possibly non-linear dynamical system to an infinite-dimensional linear system. In practice, this amounts to a simplification only when one can handle the latter numerically. Several numerical techniques have been developed in this regard. Many investigations focus on particular dynamical systems (linear systems, nonlinear systems with analytically known linearizations, ergodic systems) and their associated Koopman operator.
Numerical methods for (generalized) Fourier and Laplace analysis perform linearization of nonlinear systems close to steady states and limit cycles \cite{mezic-2013,mauroy-2013}, with a particular focus on the relation of isochrons and isostables to the eigenfunctions of the Koopman operator. These methods are useful in finding specific eigenfunctions of the dynamical system with desired properties, but are less suited for obtaining a general spectral decomposition of the Koopman operator. Giannakis \cite{giannakis-2015} describes how to estimate Koopman eigenfunctions with diffusion maps for ergodic systems.
Klus et al. \cite{klus-2017} discuss several data driven methods approximating transfer operators, including the \textit{variational approach of conformation dynamics} (VAC) and \textit{extended dynamic mode decomposition} (EDMD). Sparse identification of nonlinear dynamics (SINDy) searches for an optimal, sparse representation of the dynamics \cite{brunton-2016c,wu-2017}, requiring a large dictionary of simple building blocks.

The EDMD algorithm is an extension of Dynamic Mode Decomposition (DMD) \cite{rowley2009spectral,schmid2010dynamic}, and was developed by Williams et al. \cite{williams2014kernel,williams2015data}. The main improvement over DMD is the possibility to choose a set of observables, called a \textit{dictionary}. One can then approximate the Koopman operator as a linear map on the span of the finite set of dictionary elements. The spectral decomposition of this finite-dimensional linear map is numerically tractable and its spectral properties can approximate those of the Koopman operator. The original DMD algorithm can be interpreted as choosing only the system state as the observation of the system. By a ``careful'' choice of dictionary containing elements, beyond the system state observation functions, the EDMD algorithm is seen to have improved performance over DMD \cite{williams2015data}. 

A clear drawback of the EDMD algorithm is the need to make an {\it a priori} choice of dictionary. It is well-known that the choice significantly impacts the approximation quality of the spectral properties of the system~\cite{williams2014kernel,williams2015data,korda2017convergence}. For high-dimensional and highly non-linear systems, it is often not easy to make a judicious selection without prior information of the dynamics. 
In this paper, we aim to alleviate this issue by borrowing ideas from machine learning. We develop an iterative approximation algorithm which couples the EDMD with a {\it trainable} dictionary represented here by an artificial neural network, acting as a universal function approximator. The dictionary can be trained to adapt to the data it is presented with and this effectively reduces the need to specify a problem-dependent dictionary. 
To demonstrate the efficacy of our algorithm, we use it to perform approximate spectral decompositions of the Koopman operator for the autonomous Duffing oscillator and the Kuramoto-Sivashinsky PDE on a two-dimensional quasiperiodic and attracting invariant manifold. 

In the next section, we describe briefly the background for the Koopman operator viewpoint of dynamical systems. We also introduce the notation used throughout the paper.
Section \ref{sec:numerial method} provides a summary of the EDMD algorithm, followed by the introduction of our extension of it to adapt the dictionary elements to the data.
In section \ref{sec:applications}, we use two examples, namely the Duffing oscillator and the Kuramoto-Sivashinsky PDE, to demonstrate the efficacy of our approach.
Section \ref{sec:discussion} discusses the results with respect to accuracy and possible extensions, and section \ref{sec:conclusion} concludes the paper.

\section{Preliminaries\label{sec:preliminaries}}

%- dynamical systems, spatial isomorphism, Composition operator, Koopman operator

%The Koopman operator can be defined for discrete and continuous dynamical systems. Here, we only treat the discrete case, and use the flow map of continuous systems with a fixed time step.

% - something on the relation to Schroeders and Abels equation

\subsection{The Koopman operator}

Consider a measurable space $\mathcal{M}$ and a measure $\rho:\sigma\left(\mathcal{M}\right)\rightarrow\mathbb{R}$. 
Define a dynamical system 
\begin{equation}
x(n+1)=f\left(x(n)\right),\quad x(n)\in\mathcal{M}, \quad n\geq 0,
\label{eq:dyn_sys}
\end{equation}
on this space. The Koopman formalism focuses on the evolution of \textit{observables}, which are represented by functions on $\mathcal{M}$ in a suitable function space. Usually, we consider the Hilbert space
\[
L^2(\mathcal{M},\rho) = \{\phi:\mathcal{M}\rightarrow\mathbb{C}: \Vert\phi\Vert_{L^2(\mathcal{M},\rho)}<\infty\},
\]
where 
\[
\Vert\phi\Vert_{L^2(\mathcal{M},\rho)}=\int_{\mathcal{M}} \vert\phi(x)\vert^2 \rho(dx).
\]
where $\rho$ is a measure on $\mathcal{M}$. 

The Koopman operator, $\mathcal{K}$, is
defined as an operator acting on $L^{2}\left(\mathcal{M},\rho\right)$,
so that for an observable $\phi\in L^{2}\left(\mathcal{M},\rho\right)$,
we have
\begin{equation}
\mathcal{K}\phi=\phi\circ f.
\label{eq:koopman_dyn_sys}
\end{equation}
Intuitively, $\mathcal{K}$ describes the evolution of each observable value as driven by the dynamical system. We shall assume that $\rho$ and $\mathcal{M}$ are such that $\mathcal{K}$ is bounded. Consequently, since the Koopman operator is a bounded linear operator, it is amenable to spectral analysis. The so called {\it Koopman mode decomposition}
can be described as follows. Given a vector of observables, $O:\mathbb{\mathcal{M}\rightarrow C}^{d}$,
we can write 
\begin{equation}
O(n) = O(x(n)) = \sum_{k}\mu_{k}^n\phi_{k}\left(x(0)\right)m_{k}^{O}
\label{eq:koopman_decomp}
\end{equation}
where $\phi_{k}\in L^{2}\left(\mathcal{M},\rho\right)$ are eigenfunctions
of $\mathcal{K}$ with eigenvalues $\mu_{k}\in\mathbb{C}$. The vectors
$m_{k}^{O}\in\mathbb{C}^{d}$ are known as the Koopman modes associated
with the observable $O$ and the $k^{th}$ eigenfunction. 

In applications we are often interested in the full-state observable $O(x)=x$. Then, the Koopman mode decomposition can be viewed as a nonlinear counterpart of normal mode analysis. In this case, we interpret a finite dimensional nonlinear system \eqref{eq:dyn_sys} as a infinite dimensional linear system \eqref{eq:koopman_dyn_sys}, whose spectral evolution follows \eqref{eq:koopman_decomp}.  

\subsection{Continuous-time systems}
\label{sec:ctstime}
The Koopman formalism can be similarly applied to continuous-time systems by considering infinitesimal generators. Alternatively, we can interpret continuous-time systems as discrete-time ones by using the flow map. 
Consider the dynamical system
\[
\dot{x}(t) = g(x(t)).
\]
Let $\tau>0$ and define the flow map
\[
\Phi_\tau(x_0) := x(\tau)
\]
where $x(t)$ follows the dynamical system above with $x(0)=x_0$. Then, by defining $\tilde{x}(n):=x(n\tau)$ we obtain a discrete-time dynamical system as in \eqref{eq:dyn_sys}, with $f\equiv\Phi_\tau$. 

\section{Numerical methods}\label{sec:numerial method}
Although the linear Koopman dynamics is theoretically easier to analyze, in practice it is often a challenge to compute its spectral properties due to its infinite-dimensionality. In this section, we briefly review the classical EDMD method for computing Koopman mode decompositions. We then introduce our numerical method that incorporates machine learning to address the important shortcoming of traditional methods - the need of selecting a fixed dictionary that may be generally ``inefficient''. 

\subsection{The EDMD algorithm}
Since our dictionary learning algorithm is built upon the extended dynamic mode decomposition (EDMD) \cite{williams2015data}, we begin by describing briefly the EDMD algorithm. The main idea is to estimate a finite-dimensional representation of the Koopman operator $\mathcal{K}$ in the form of a finite-dimensional linear map $K$, whose spectral properties will then approximate those of $\mathcal{K}$. To do this, pick a dictionary $\Psi = \{\psi_{1}, \psi_{2}, \dots, \psi_{M}\}$ where
\[
\psi_{i}: \mathcal{M} \rightarrow \mathbb{R} \quad \text{for} \quad i=1,2,\dots,M.
\]
Now, consider the span $U(\Psi)=\text{span}\{\psi_1,\dots,\psi_M\}=\left\{a^T \Psi :a\in\mathbb{C}^M \right\}$, which is a linear subspace of $L^2(\mathcal{M},\rho)$. For any $\phi= a^T \Psi \in U(\Psi)$, we have 
$\mathcal{K}\phi = a^T \mathcal{K}\Psi = a^T \Psi\circ f$.
If $\mathcal{K}(U(\Psi))=U(\Psi)$, then we also have $\mathcal{K}\phi = b^T \Psi$ for some $b\in\mathcal{C}^M$. Hence, a finite dimensional representation of $\mathcal{K}$ is realized as the matrix $K\in\mathbb{R}^{M\times M}$ with $b = K^T a$. Thus, we have the equality $a^T\Psi\circ f=a^T K\Psi$. For this to hold for all $a$ we must have $\Psi\circ f = K\Psi$. To find $K$, we use pairs of data points $\left\{ x(n),y(n)\right\} _{n=1}^{N}$ with $y(n)=f\left(x(n)\right)$ and solve the minimization problem
\begin{equation}
K = \operatornamewithlimits{arg min}_{\tilde{K}\in\mathbb{R}^{M\times M}} J(\tilde{K}) = \sum_{n=1}^{N} {\Vert \Psi(y(n)) - \tilde{K}\Psi(x(n)) \Vert^2}.
\label{eq:edmd_min_prob}
\end{equation}
If $U(\Psi)$ is indeed invariant under $\mathcal{K}$, then $J(K)=0$. Otherwise, $J(K)>0$ and the procedure above seeks to find $K$ that minimizes the residual $J$. The solution to \eqref{eq:edmd_min_prob} is
\[
K = G^{+} A
\]
with
\begin{align}
	G &= \frac{1}{N}\sum_{n=1}^{N}\Psi\left(x(n)\right)^{T}\Psi\left(x(n)\right), \nonumber\\
	A &= \frac{1}{N}\sum_{n=1}^{N}\Psi\left(x(n)\right)^{T}\Psi\left(y(n)\right), 
	\label{eq:GA}
\end{align}
and $G^{+}$ denotes the pseudo-inverse of $G$.
With $K$ derived, it is straightforward to find eigenfunctions and eigenvalues of $K$, and likewise modes associated with an observable $O$. 
For example, one can check that for each right eigenvector $\xi_{j}$ of $K$ with eigenvalue $\mu_{j}$, the function
\[
\phi_{j}=\xi_{j}^{T}\Psi
\]
is an approximation of an eigenfunction of $\mathcal{K}$ with the same eigenvalue $\mu_j$. Also, for any vector of observables $O=B\Psi$, the $j^{th}$ Koopman mode associated with $u$ is given by 
\[
m_j = B\zeta_j
\]
where $\zeta_j$ is the $j^{th}$ left eigenvector of $K$. 

The matrix $K$ found in this way is shown to converge to $\mathcal{K}_\Psi$, the $L^2$ orthogonal projection of the Koopman operator $\mathcal{K}$ onto $U(\Psi)$, as $N\rightarrow\infty$ \cite{rowleyweb17,williams2015data}. It has been further established that if $\Psi$ consists of linearly independent basis functions, then as $M\rightarrow\infty$ one has $\mathcal{K}_\Psi \rightarrow \mathcal{K}$ in the strong operator topology \cite{korda2017convergence}. 

In practice, however, both $N,M$ are finite. Therefore, one primary assumption underlying EDMD's application is that the finite dimensional subspace $U(\Psi)$ is approximately invariant under $\mathcal{K}$. This is true if either $M$ is very large or more practically, if the dictionary set $\Psi$ is judiciously chosen \cite{williams2015data}. The choice of dictionary is especially difficult for highly nonlinear or high dimensional systems, for which even enumerating a standard basis (e.g. orthogonal polynomials) becomes prohibitively expensive. Although there exist kernel methods \cite{williams2014kernel} to alleviate such problems, the choice of dictionary sets (including kernel functions) remains a central challenge for the general applicability of EDMD. 

In next section, we use ideas from machine learning to show how one can alleviate the problem of having to choose a fixed dictionary. Most importantly, this holds the promise of allowing high-quality representation with relatively fewer dictionary terms. 

\subsection{EDMD with Dictionary Learning (EDMD-DL)}

Dictionary learning (or sparse coding) is a classical problem in signal
processing and machine learning \cite{donoho2006compressed}. The problem statement is as follows: given a set of input data,
$X=\left(x(1)\quad x(2)\quad...\quad x(N)\right)\in\mathbb{R}^{d\times N}$, 
we wish to find a sparse representation of it in the form of $X=DK$, where $D\in\mathbb{R}^{d\times M}$
is a size-$M$ set of dictionary vectors and $K\in\mathbb{R}^{M\times N}$ is a sparse representation. 
For any fixed $D$, it is difficult to fulfill accuracy ($X\approx DK$) {\it and} sparsity ($\Vert K\Vert_0$ small)
at the same time. A better approach is to make $D$ adapted to data and solve
\begin{equation}
(K,D) = \operatornamewithlimits{arg min}_{(\tilde{K},\tilde{D})} J(\tilde{K},\tilde{D}) =  \Vert X-\tilde{D}\tilde{K}\Vert_{F}^{2}+\lambda\Vert \tilde{K}\Vert_{1}
\label{eq:var_dict_learning}
\end{equation}
where $\Vert\cdot\Vert_{F}$ is the Frobenius norm. The $\ell_{1}$
penalty induces sparsity without turning it into a combinatorial optimization problem, as in the case for $\ell_0$ penalty. To remove degeneracies, one may impose further conditions such as $\left\{ D:\Vert D\Vert_{F}=1\right\}$. 

From the above, the key idea we would like to adapt to the EDMD framework is allowing the dictionary to be trainable, which then enables one to find an ``efficient'' representation of the data with a smaller number of adaptive basis elements. In this sense, the procedure is similar to the Karhunen-Lo\`{e}ve decomposition (KLD) \cite{karhunen1947lineare,loeve1978probability}, whose sampled versions are also known as principal component analysis (PCA) and proper orthogonal decomposition (POD). The goal of KLD is to obtain an expansion of stochastic processes in terms of adaptive basis functions for which the truncation error is optimal, in the mean-squared sense. 

In our case of EDMD decompositions, our goal is to make the dictionary set $\Psi$ adaptive so that we can minimize the norm of the residual $\Psi\circ f-K\Psi$ resulting from the finite-dimensional projection (see \eqref{eq:edmd_min_prob}). Hence, instead of \eqref{eq:edmd_min_prob} we can consider the extended minimization problem
\begin{align}
	(K, \Psi) =& \operatornamewithlimits{arg min}_{(\tilde{K}, \tilde{\Psi})} J(\tilde{K},\tilde{\Psi}) \nonumber\\
	=& \sum_{n=1}^{N} {\Vert \tilde{\Psi}(y(n)) - \tilde{K}\tilde{\Psi}(x(n)) \Vert^2} + \lambda(\tilde{K},\tilde{\Psi}),
	\label{eq:dl_gen_min_prob}
\end{align}
where $\lambda(K,\Psi)$ is a suitable regularizer. Unlike \eqref{eq:var_dict_learning}, the dictionary functions in $\Psi$ are not assumed to be linear functions and hence nonlinear optimization methods must be used. Nevertheless, provided we can solve \eqref{eq:dl_gen_min_prob}, this formulation provides us with a means to find an adaptive set of dictionary elements that give optimal truncation errors, in a similar vein to sparse coding or the Karhunen-Lo\`{e}ve decomposition.

We have outlined the primary idea underlying our adaptive EDMD algorithm. Next, we present a computational algorithm to solve \eqref{eq:dl_gen_min_prob} by parameterizing it with neural networks. 

%There are many possible solution methods for the problem \eqref{eq:var_dict_learning} \cite{engan1999method}, \cite{lee2007efficient}.
%The simplest of which that is applicable to very large scale problems
%is (stochastic) gradient descent \cite{ruder2016overview}, which consists of solving \eqref{eq:var_dict_learning} by alternating between two steps: first we optimize $D$ holding $K$ constant, and then we optimize $K$ holding $D$ constant. This is iterated until convergence. This idea forms the basis of our extension to the EDMD algorithm, which we now outlined in the next section.  

\subsection{A practical algorithm}
\begin{figure}
	\includegraphics[width=0.8\textwidth]{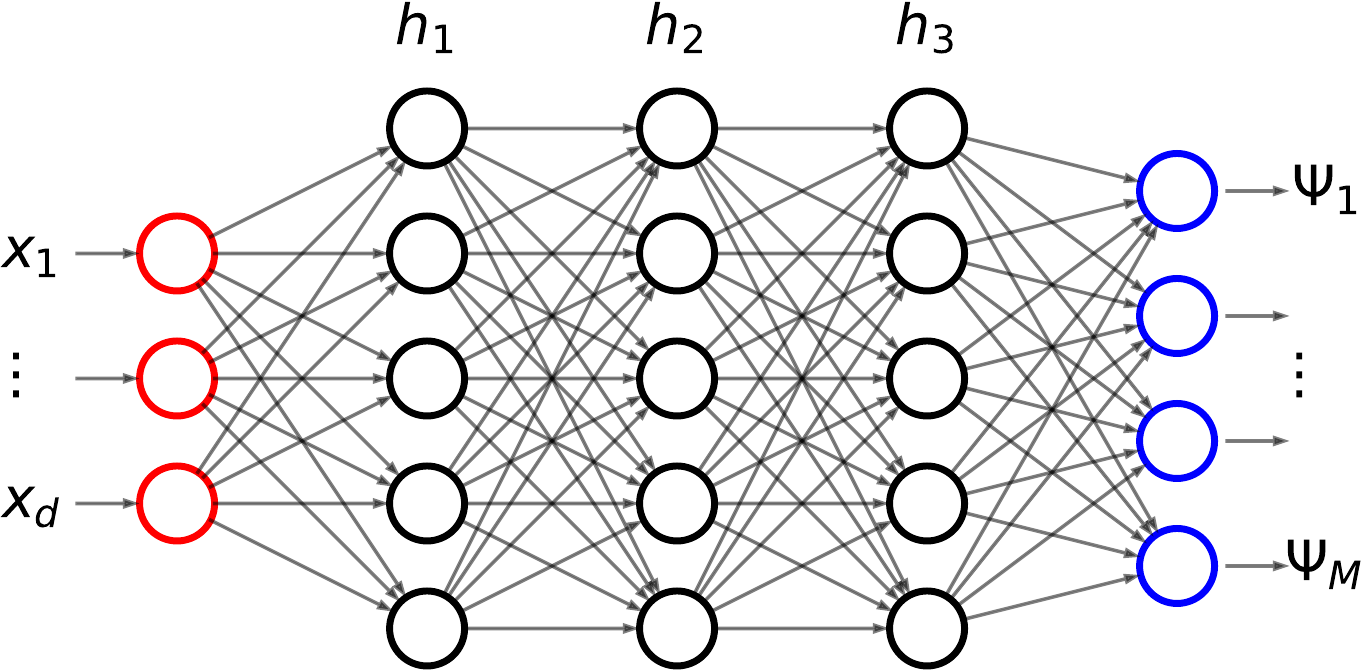}
	\caption{Neural network function approximator for the trainable dictionary $\Psi(x;\theta)$. The network is fully connected and consists of 3 hidden layers $h_1,h_2,h_3$. Arrows connecting layers corresponds to affine transformations followed by $\tanh$ activations. See Eq. \eqref{eq:nn}.}
	\label{fig:nn}
\end{figure}

To solve \eqref{eq:dl_gen_min_prob}, we parameterize $\Psi$ by a universal function approximator, i.e. $\Psi(x) = \Psi(x;\theta)$ for some $\theta\in \Theta$ to be varied. Here, we choose a simple feed-forward, 3-layer neural network as the approximator for $\Psi$ (see Fig. \ref{fig:nn}). Concretely, we choose a hidden dimension $\ell$ and set
\begin{align}
	\Psi(x) &= W_{\text{out}}h_3 + b_{\text{out}}, \nonumber\\
	h_{k+1} &= \tanh(W_{k}h_{k} + b_k), \quad k=0,1,2
	\label{eq:nn}
\end{align}
where $h_0 = x$ and $W_{0}\in\mathbb{R}^{\ell\times d}, b_{0}\in\mathbb{R}^{d}$, $W_{\text{out}}\in\mathbb{R}^{M\times \ell}, b_{\text{out}}\in\mathbb{R}^M$ and $W_{k}\in\mathbb{R}^{\ell\times \ell}, b_{k}\in\mathbb{R}^{\ell}$ for $k=1,2$. The set of all trainable parameters is $\theta = \{W_{\text{out}},b_{\text{out}},\{W_k,b_k\}_{k=0}^2\}$, which contains a total of $d (l+1)+l (2 l+M+3)$ scalar variables.

With $\Psi$ parameterized, we can then solve \eqref{eq:dl_gen_min_prob}: 
\begin{align}
	(K, \theta) =& \operatornamewithlimits{arg min}_{(\tilde{K}, \tilde{\theta})} J(\tilde{K},\tilde{\theta}) \nonumber\\
	=& \sum_{n=1}^{N} {\Vert \Psi(y(n);\tilde{\theta}) - \tilde{K}\Psi(x(n);\tilde{\theta}) \Vert^2} + \lambda \Vert \tilde{K}\Vert^2_F.
	\label{eq:dl_min_prob}
\end{align}
We picked the Tikhonov regularization \cite{tikhonov1943stability,ng2004feature} with identity matrix for $\tilde{K}$ to improve the stability of the algorithm. 
Notice that if there exists $\tilde{\theta}$ with $\Psi(x;\tilde{\theta})\equiv0$, the right hand
side identically vanishes and the minimum is trivially attained. 
Thus, to obtain meaningful approximations we need further restrictions.
A natural one is to include in $\Psi=\{\psi_1,\dots,\psi_M\}$
some fixed (non-trainable) functions, such as the constant and the projection maps. The presence of the latter is important for reconstructing trajectories. This is because to find the Koopman modes we require the identity map $O(x)=x$, whose components are projection maps, to be in the linear span $U(\Psi)$. The inclusion of these non-trainable dictionary functions then removes the possibility that 
$\Psi(x;\tilde{\theta})\equiv0$. 

We solve \eqref{eq:dl_min_prob} by iterating the following two steps: (a) Fix $\theta$, optimize $K$; Then (b) fix $K$, optimize $\theta$.

\begin{enumerate}
	\item[(a)] \textbf{Fix $\theta$, optimize $K$.}
	For fixed $\theta$ (hence fixed $\Psi$), \eqref{eq:dl_min_prob} is almost the same problem as \eqref{eq:edmd_min_prob}, but with the addition of the Tikhonov regularizer. The solution is \cite{golub2012matrix}
	\begin{equation}
	\tilde{K} = (G(\tilde{\theta})+\lambda I)^{+} A(\tilde{\theta})
	\end{equation}
	where $G,A$ are defined in \eqref{eq:GA} and $I$ is the $d$-dimensional identity matrix. 
	\item[(b)] \textbf{Fix $K$, optimize $\theta$.} 
	This is a standard machine learning problem. As there is no linear
	structure in the problem, we cannot write down its exact solution.
	Instead, we proceed by iterative updates in the form of gradient descent, i.e., we set 
	\begin{equation}
	\tilde{\theta} \leftarrow \tilde{\theta} - \delta \nabla_{\theta} J(\tilde{K}, \tilde{\theta}).
	\end{equation}
	If both the dimension $d$ and the sample size $N$ is large, $\nabla_\theta J$ will be expensive to evaluate.
	We can then employ stochastic gradient descent and its variants \cite{ruder2016overview}, where the gradient
	$\nabla_\theta J$ is replaced by randomly sampled unbiased estimators. 
\end{enumerate}
The above two steps are iterated until convergence. We have observed (empirically!) that the algorithm performed stably and converged for general initializations. A rigorous proof of the convergence of the algorithm will be left as future work. 
The algorithm is summarized in Alg. \ref{alg:EDMDdict} and we hereafter refer to it as EDMD with dictionary learning (EDMD-DL). 
\begin{algorithm}[H]
	Initialize $K,\theta$. 
	
	Set learning rate $\delta>0$, tolerance $\epsilon>0$, regularizer $0<\lambda\ll 1$
	
	\textbf{while} $J(K,\theta)>\epsilon$ \textbf{do}:
	\begin{itemize}
		\item[] $K \leftarrow  (G(\theta)+\lambda I)^{-1} A(\theta) $ 
		\item[] $\theta \leftarrow \theta - \delta \nabla_\theta J(K, \theta)$
	\end{itemize}
	\caption{EDMD with dictionary learning (EDMD-DL) \label{alg:EDMDdict}}
\end{algorithm}

\section{Applications of EDMD-DL\label{sec:applications}}

In this section, we compare the results from the EDMD-DL algorithm with the classical EDMD results on various example problems to 
illustrate the advantages of an adaptive, trainable dictionary. For each example,
we evaluate the performance of various methods by two quantitative metrics:
\begin{itemize}
	\item \textbf{Accuracy of trajectory reconstruction}. We reconstruct trajectories using the Koopman mode decomposition formula \eqref{eq:koopman_decomp} with $O(x)=x$. We then monitor the reconstruction error as
	\begin{equation}
	\text{Error} = \sqrt{\frac{1}{N}\sum_{n=1}^{N} \vert x(n) - \tilde{x}(n) \vert^2},
	\label{eq:recon_error_def}
	\end{equation}
	where $x$ is the true trajectory data (according to \eqref{eq:dyn_sys}) and $\tilde{x}$ is the reconstructed trajectory (according to \eqref{eq:koopman_decomp}). 
	\item \textbf{Accuracy of eigenfunction approximation}. For each $j=1,2,\dots,M$ we define the eigenfunction approximation error
	\begin{equation}
	E_j = \Vert \phi_j\circ f - \mu_j \phi_j \Vert_{L^2(\mathcal{M},\rho)},
	\label{eq:Ej}
	\end{equation}
	where $\phi_j$ and $\mu_j$ are the $j^\text{th}$ eigenfunction and eigenvalue found by the algorithm. The quantity above can be approximated by Monte-Carlo integration
	\[
	E_j \approx \sqrt{\frac{1}{I}\sum_{i=1}^I \vert \phi_j\circ f(x(i)) - \mu_j \phi_j(x(i)) \vert^2},
	\]
	where $x(i)\sim \rho$ are identically and independently distributed for all $i$. 
\end{itemize}

%- briefly discuss the error metrics
% * error by comparing predicted trajectories to real ones
% * compare error against number of observables
% * error for individual eigenfunctions

\subsection{Duffing oscillator}

We start by applying EDMD-DL to the Duffing oscillator, which describes the evolution of $x=(x_1, x_2)$ governed by
\begin{align}
	\dot{x_1} &= x_2, \nonumber \\
	\dot{x_2} &= -\delta x_2 - x_1 (\beta + \alpha x_1^2). 
	\label{eq:duffing}
\end{align}
We take $\alpha=1$, $\beta=-1$ and $\delta=0.5$ so that there are two stable steady states at $(\pm 1, 0)$ separated by a saddle point at $(0,0)$. We convert the continuous dynamical system to a discrete one by defining flow maps as discussed in \ref{sec:ctstime}, with the choice $\tau=0.25$. We draw 1000 random initial conditions uniformly in the region $[-2,2]^2$. Each initial condition is evolved up to $n=10$ steps with the flow-map so that we have a total of $10^5$ data points to form the training set. 

\begin{figure}
	\begin{center}
		\includegraphics[width=0.8\textwidth]{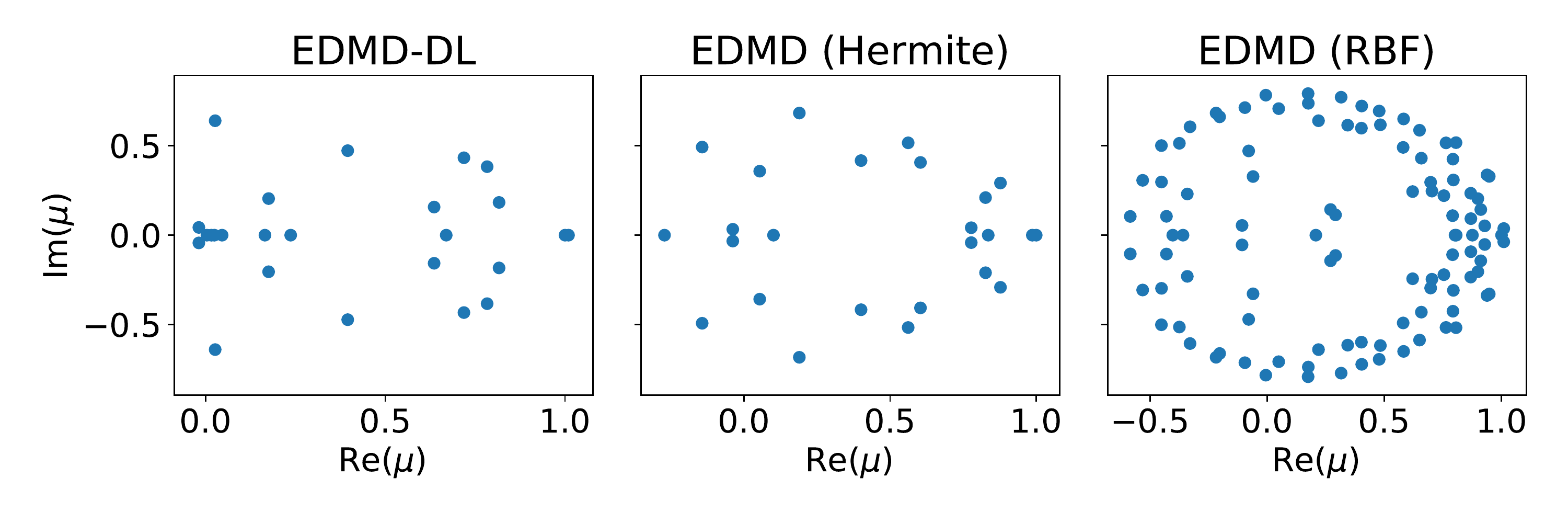}
	\end{center}
	\caption{Eigenvalues of the Koopman operator for the Duffing oscillator estimated from each algorithm. For EDMD-DL and EDMD with Hermite basis, 25 estimated eigenvalues are shown (since 25 dictionary functions are used). For EDMD with RBF dictionary, 100 RBF functions are used and so 100 estimated eigenvalues are shown. We see that EDMD with Hermite polynomials has found many more eigenvalues with large magnitudes but does not improve accuracy significantly over EDMD-DL (see Fig. \ref{fig:duffing_traj}). In this sense, we see that EDMD-DL has found a more efficient representation. }
	\label{fig:duffing_evals}
\end{figure}

\begin{figure}
	\begin{center}
		\includegraphics[width=0.8\textwidth]{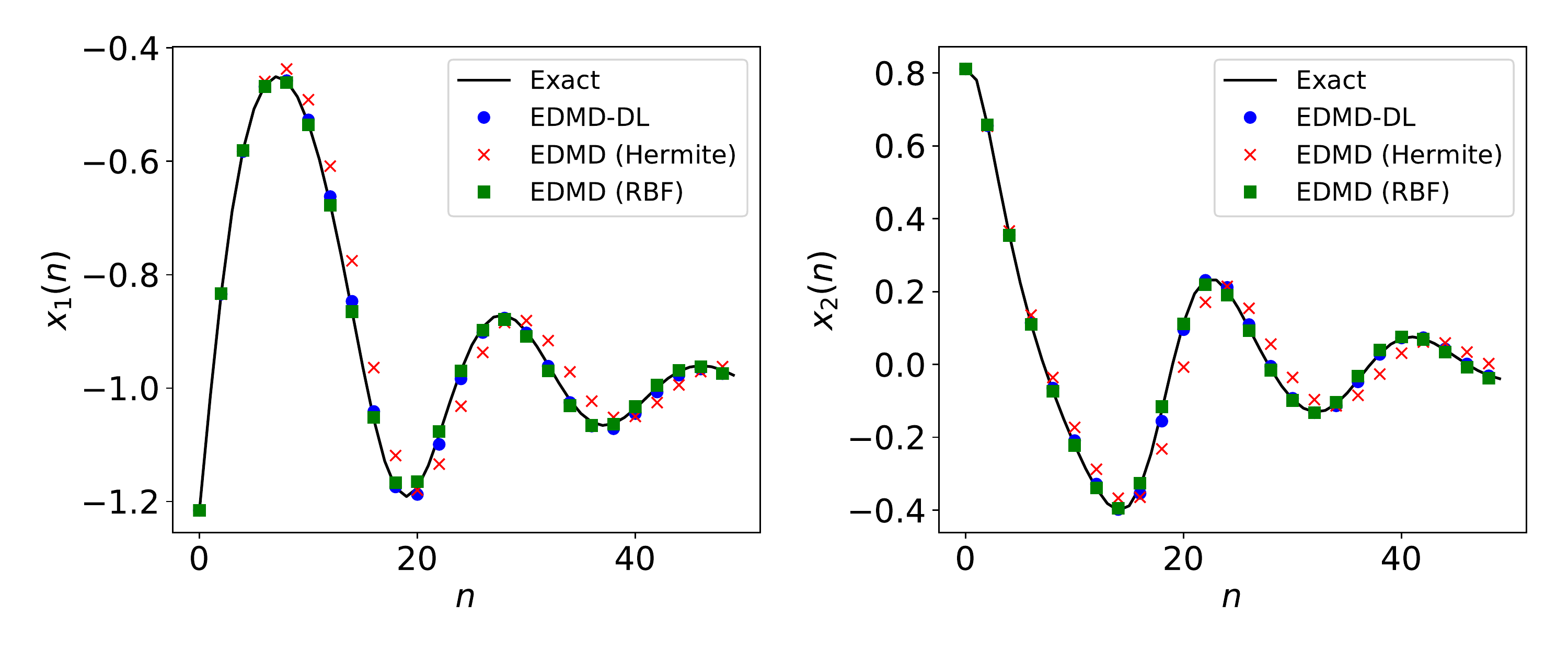}
		
		\includegraphics[width=0.8\textwidth]{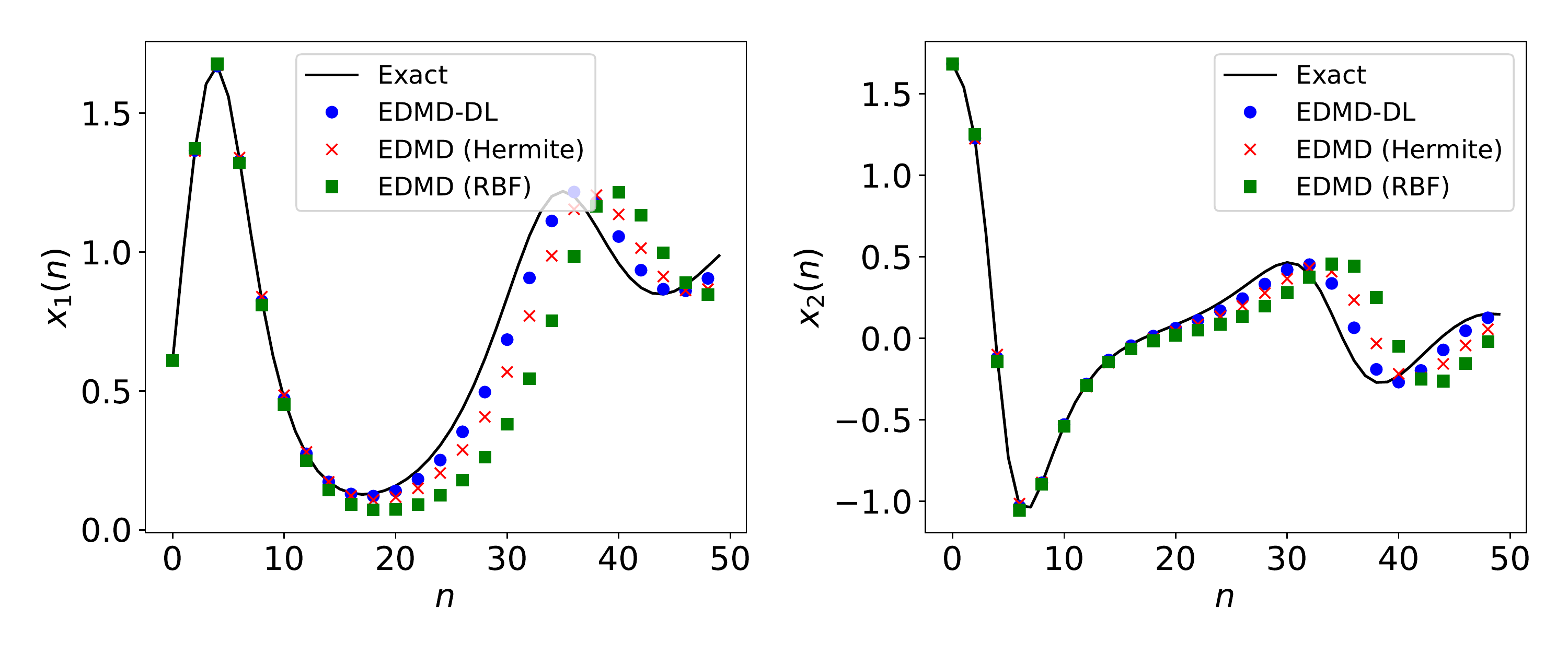}
		
		\includegraphics[width=0.8\textwidth]{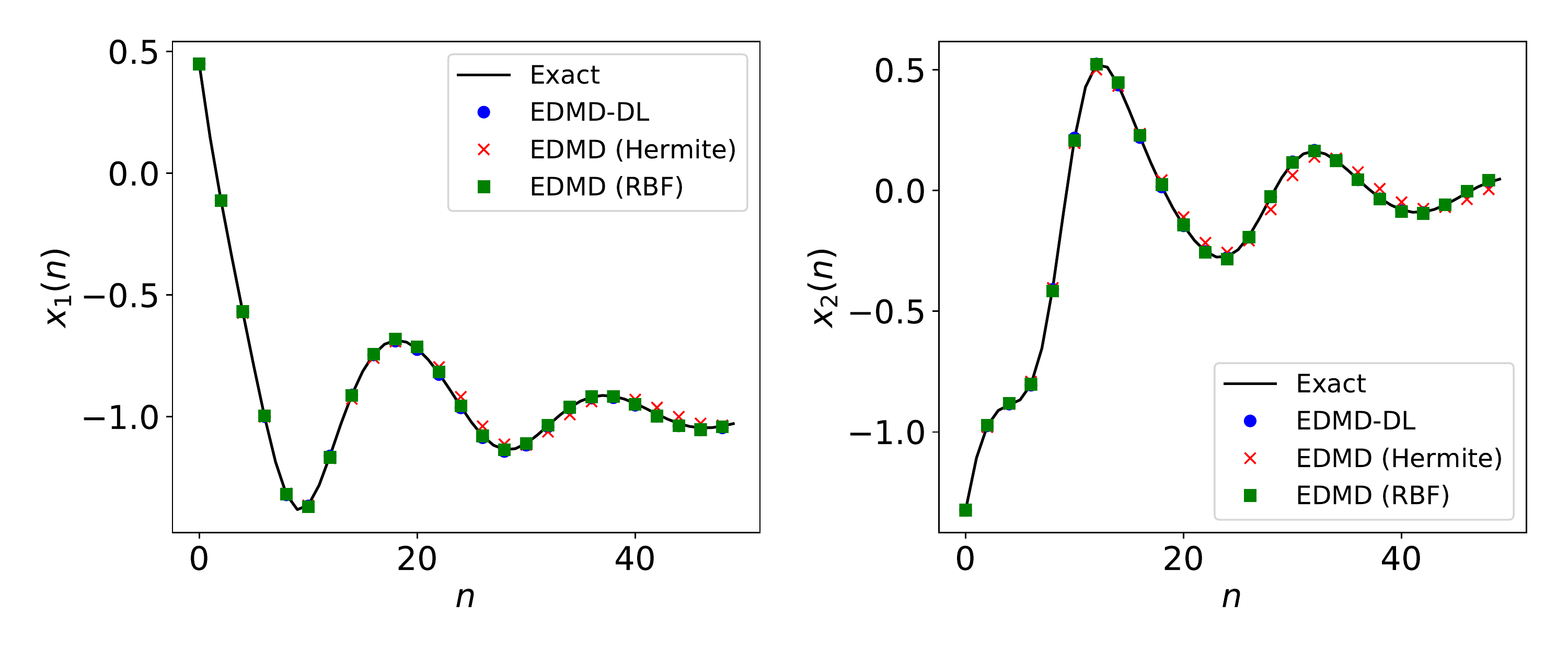}
	\end{center}
	\caption{Trajectories of the Duffing oscillator reconstructed from Koopman decomposition using various algorithms. Three different initial conditions in $[-2,2]^2$ are selected. We observe that EDMD-DL (with 25 dictionary elements) has better reconstruction accuracy than classical EDMD with Hermite polynomials with the same number of dictionary elements. We also see that EDMD-DL performs approximately on par with EDMD with RBF dictionary (100 dictionary elements), which is known to be especially suited for this problem \cite{williams2015data}. A quantitative comparison of reconstruction errors vs dictionary size is given in Fig. \ref{fig:num_dict}(a)}
	\label{fig:duffing_traj}
\end{figure}

\begin{figure}
	\begin{center}
		\subfloat[Duffing]{\includegraphics[width=0.4\textwidth]{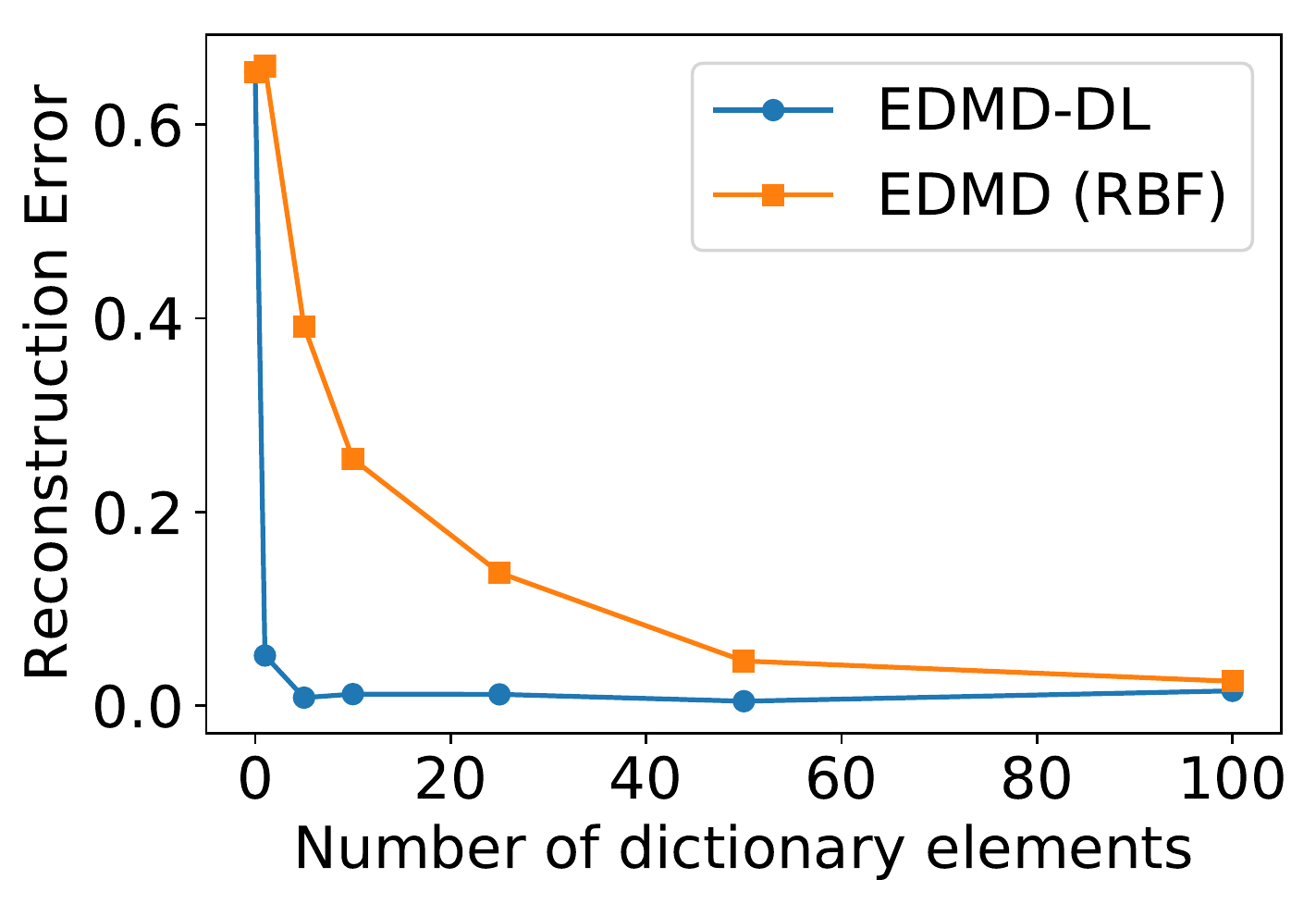}}
		\subfloat[Kuramoto-Sivashinsky]{\includegraphics[width=0.4\textwidth]{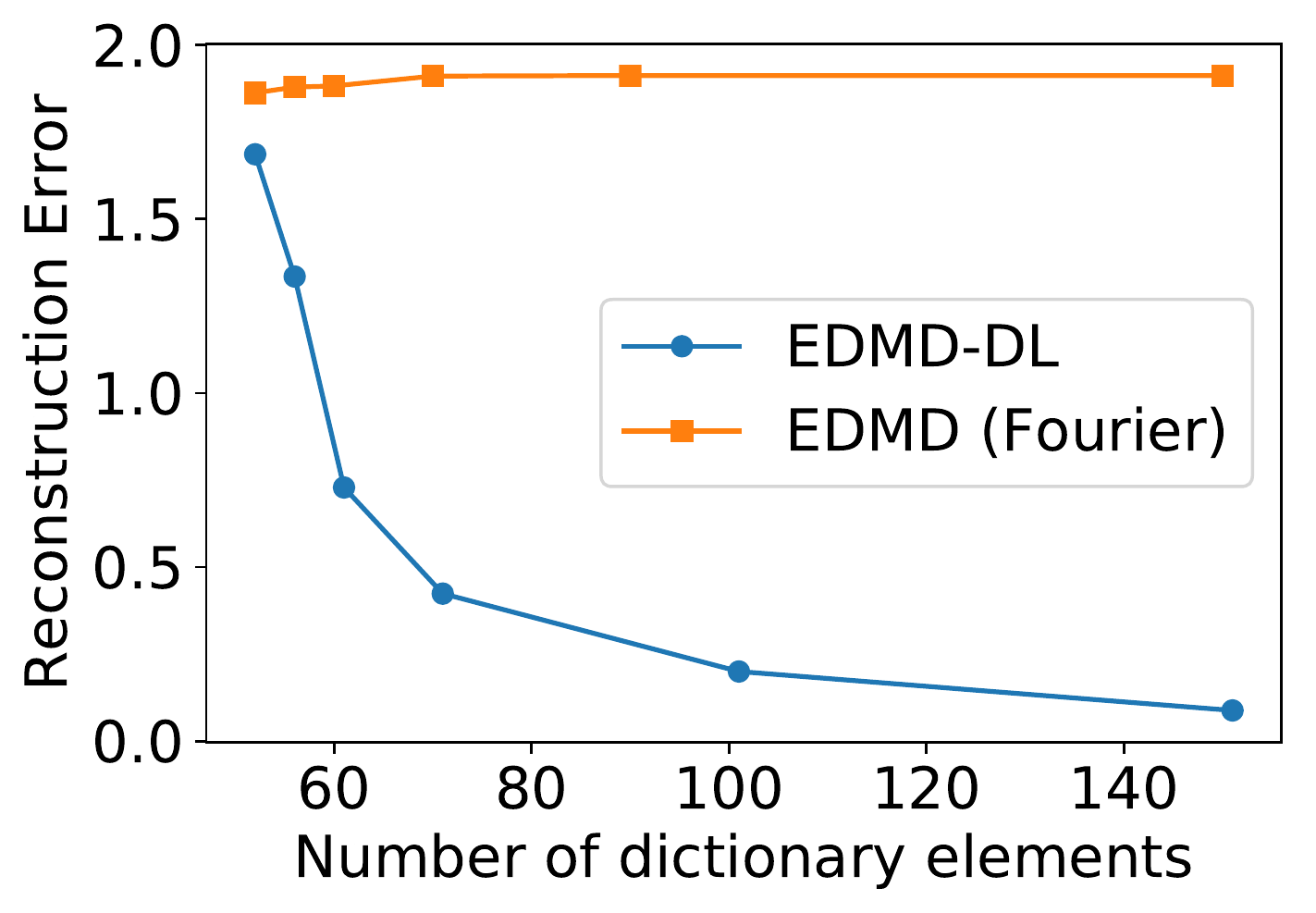}}
	\end{center}
	\caption{Trajectory reconstruction error for EDMD-DL vs classical EDMD with hand-picked dictionary, whose sizes are varied. The errors are averaged over 50(10) random and unseen initial conditions for the Duffing (Kuramoto-Sivashinsky) system. We see that EDMD-DL requires much smaller dictionary sizes in order to capture the system's dynamics. }
	\label{fig:num_dict}
\end{figure}

Now, we apply the EDMD-DL algorithm with 22 trainable dictionary outputs (plus 3 non-trainable ones, i.e. one constant map and two coordinate projection maps) and compare its performance to EDMD with two choices of dictionary sets 1) using 25, two-dimensional Hermite polynomials, and 2) 100 thin-plate radial basis functions (RBF) with centers placed on the training data using k-means clustering (\texttt{scipy.cluster.vq.kmeans}, thin-plates $r^2\ln(r+\delta)$, regularized with $\delta=10^{-4}$). In Fig. \ref{fig:duffing_evals}, we show the eigenvalues found by the three methods. To quantitatively compare the performance, we plot the trajectories reconstructed by the Koopman decomposition against the exact trajectories obtained by integrating the evolution equations \eqref{eq:duffing}. The results are shown in Fig. \ref{fig:duffing_traj}. We see that although EDMD-DL uses a small set of trainable dictionary outputs, it out-performs both EDMD with Hermite polynomials and RBFs, despite the fact that the latter is carefully chosen to be effective for the Koopman decomposition of the Duffing equation \cite{williams2015data}. To confirm that EDMD-DL requires a smaller dictionary size, we plot in Fig. \ref{fig:num_dict}(a) the reconstruction error averaged over 50 random initial conditions vs the dictionary size for EDMD-DL and EDMD with RBF dictionaries. We see that EDMD-DL achieves lower reconstruction error at smaller dictionary sizes. 

\begin{figure}
	\begin{center}
		\includegraphics[width=0.8\textwidth]{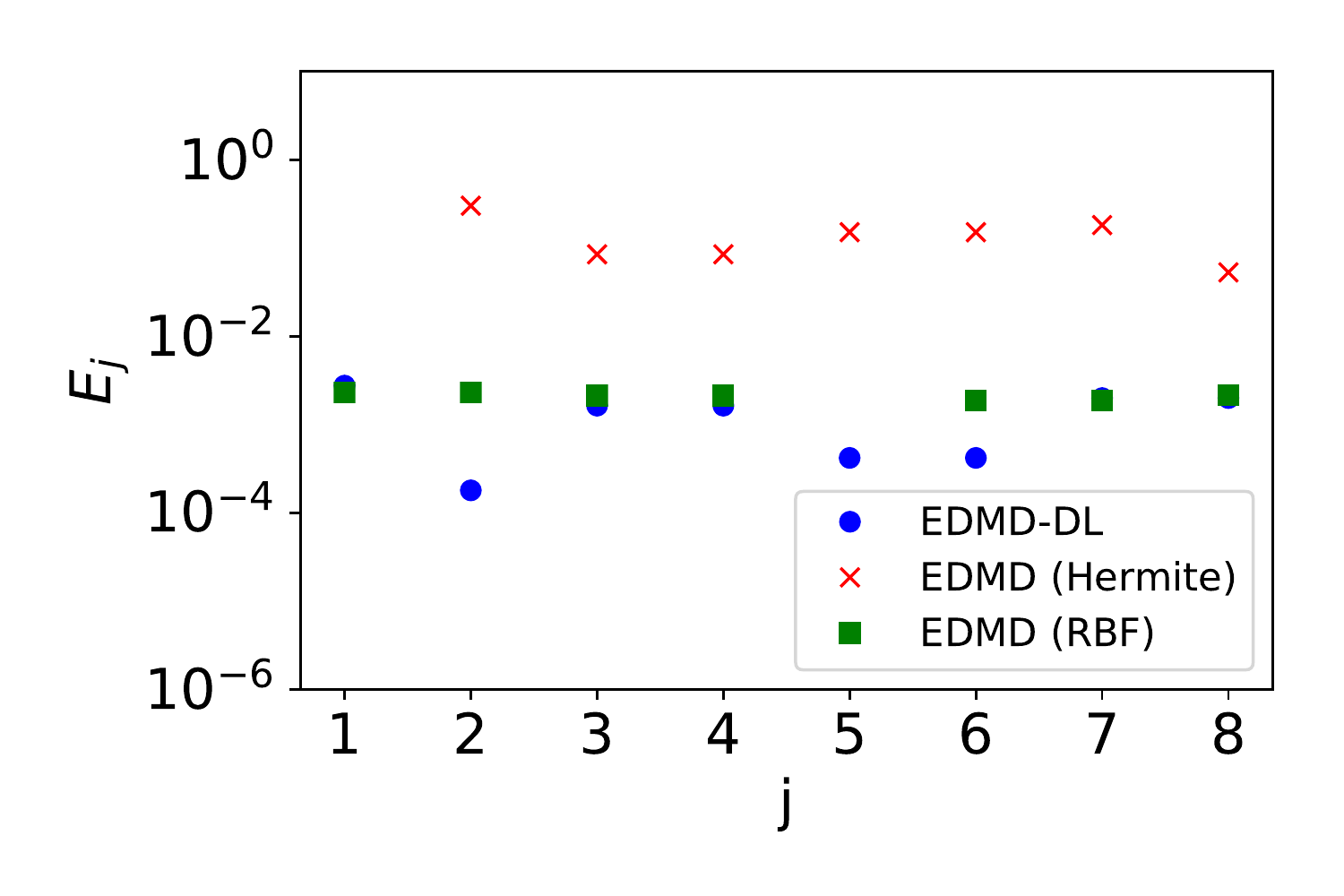}
	\end{center}
	\caption{Eigenfunction errors for the Duffing oscillator. Both EDMD-DL and EDMD with Hermite dictionary have 25 dictionary functions. EDMD with RBF dictionary has 100 dictionary elements. Again, we observe that dictionary learning has comparable performance to the well-chosen (and large) RBF dictionary and has better performance than the Hermite dictionary. }
	\label{fig:duffing_efuns}
\end{figure}

As a further quantitative comparison, we evaluate the quality of the eigenfunctions by calculating for each $j$ the eigenfunction error $E_j$ (See definition \eqref{eq:Ej}) with $\rho = 1_{[-2,2]^2}$. The value of $E_j$ for the first 8 leading eigenfunctions are shown in Fig. \ref{fig:duffing_efuns}. Again, we can see that EDMD-DL achieves comparable performance with a well-picked dictionary (RBF), and outperforms poorly picked ones (Hermite). 

The Duffing oscillator is a low dimensional dynamical system, hence enumerating polynomial basis functions is still reasonably tractable.  Provided that enough of them are included in the dictionary, the finite-dimensional approximations for the Koopman operator become reasonably accurate. Moreover, {\it a priori} domain knowledge of the eigenfunctions can also allow us to pick better dictionaries, such as the RBFs \cite{williams2015data}. Consequently, standard EDMD is still reasonable even though EDMD-DL still performs better. 
For general high dimensional systems, it is difficult to choose a dictionary in a systematic and efficient way. This situation is precisely where dictionary learning is most advantageous. 

\subsection{Kuramoto-Sivashinsky PDE}
Consider the Kuramoto-Sivashinsky PDE
\begin{equation}
u_t + 4 u_{zzzz} + \alpha( u_{zz} + uu_z ) = 0, \quad z\in[0, 2\pi]
\end{equation}
with $\alpha=16$ and periodic boundary conditions on $\left[0,2\pi\right)$. The initial condition is parameterized with $a\in[0.8,1]$, $b\in[0.5,1]$, and given as
\begin{equation}
u(z,0)=a\sin(2\pi z)+b\exp(\cos(2\pi z)).
\end{equation}
We sample $a$ and $b$ with a 100 random, uniformly distributed points, and compute the solution $u(z,t)$ at 50 equally distributed spatial points on $\left[0,2\pi\right)$, and 100 points in time, in the interval $\left[0,0.5\right]$.

\begin{figure}
	\begin{center}
		\includegraphics[width=0.8\textwidth]{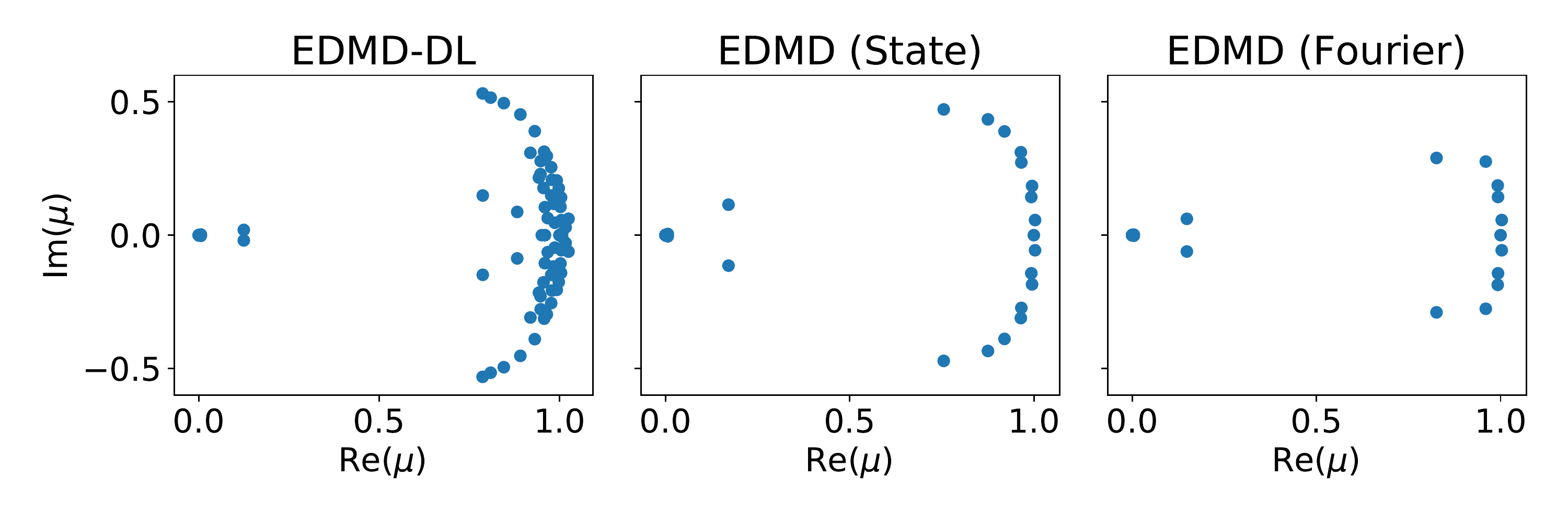}
	\end{center}
	\caption{Eigenvalues of the Koopman operator of the Kuramoto-Sivashinsky PDE estimated from each algorithm. Number of eigenvalues correspond to dictionary sizes, which are: 101 for EDMD-DL; 250 EDMD with dictionary containing states and derivatives; 150 for EDMD with dictionary containing Fourier coefficients. Observe that although EDMD-DL produced less eigenvalues (because of a smaller dictionary), it produced more meaningful eigenvalues as compared to those of EDMD, where most are concentrated at 0. This is the opposite case to Fig. \ref{fig:duffing_evals} because the PDE system necessarily requires a richer representation. Although both EDMD (state) and EDMD (Fourier) produced less eigenvalues with large magnitudes, they result in inaccurate representations of the dynamics (see Fig. \ref{fig:kse_traj}) and hence cannot be considered sparse representations. }
	\label{fig:kse_evals}
\end{figure}

\begin{figure}
	\begin{center}
		\includegraphics[width=0.8\textwidth]{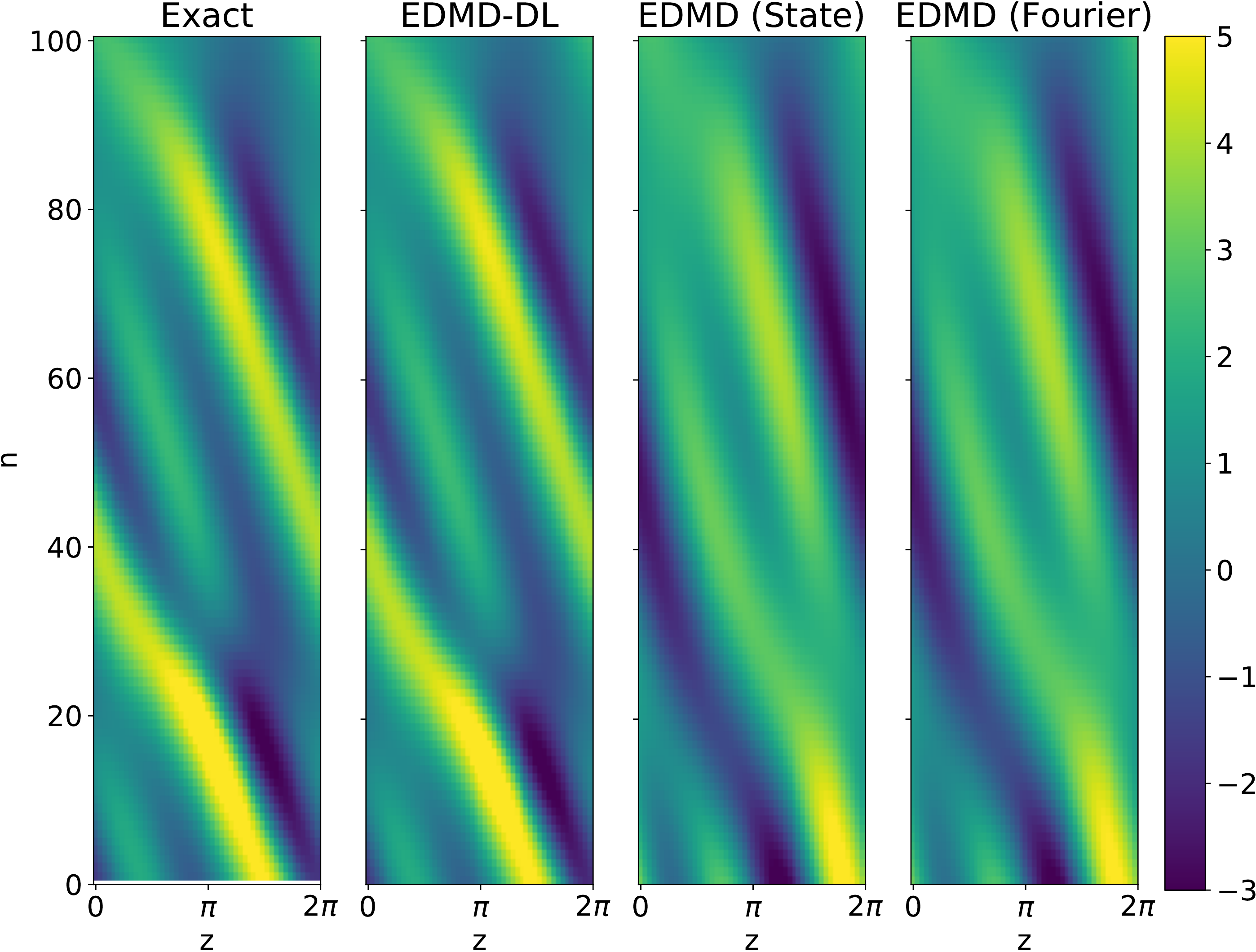}
	\end{center}
	\caption{Trajectory of the Kuramoto-Sivashinsky PDE reconstructed from Koopman decomposition using various algorithms. A random initial condition is picked from the same distribution as, but not from, the training data. We see that classical EDMD with both dictionaries cannot reproduce fine structures of the solution, whereas EDMD-DL performs well by adapting the dictionary to data. Also see Fig. \ref{fig:num_dict} for a quantitative comparison of the reconstruction error vs number of dictionary functions. }
	\label{fig:kse_traj}
\end{figure}

As discussed before, it is difficult to pick a dictionary for the classical EDMD algorithm. Here, we use two choices:
\begin{enumerate}
	\item A dictionary containing the state $u(z_i,t)$ and four of its spatial derivatives, all sampled at 50 equally spaced grid points $z_i\in\left[0,2\pi\right)$. Thus, in total, this dictionary contains 250 elements.
	\item A dictionary containing the state $u(z_i,t)$, sampled at 100 equally spaced grid points $z_i\in\left[0,2\pi\right)$, and its Fourier coefficients, separated into real and imaginary parts. In total, this dictionary contains 150 elements.
\end{enumerate}
These two dictionaries are compared to the results of EDMD-DL, where we pick 50 trainable dictionary outputs on top of the constant and projection maps (so that $M=101$). In Fig. \ref{fig:kse_evals}, we plot the eigenvalues found by each algorithm. We observe that although EDMD-DL uses a smaller number of dictionary outputs, the eigenvalue spectrum found is richer than those found by classical EDMD, where many computed eigenvalues are effectively $0$. 
\begin{figure}
	\begin{center}
		\includegraphics[width=0.8\textwidth]{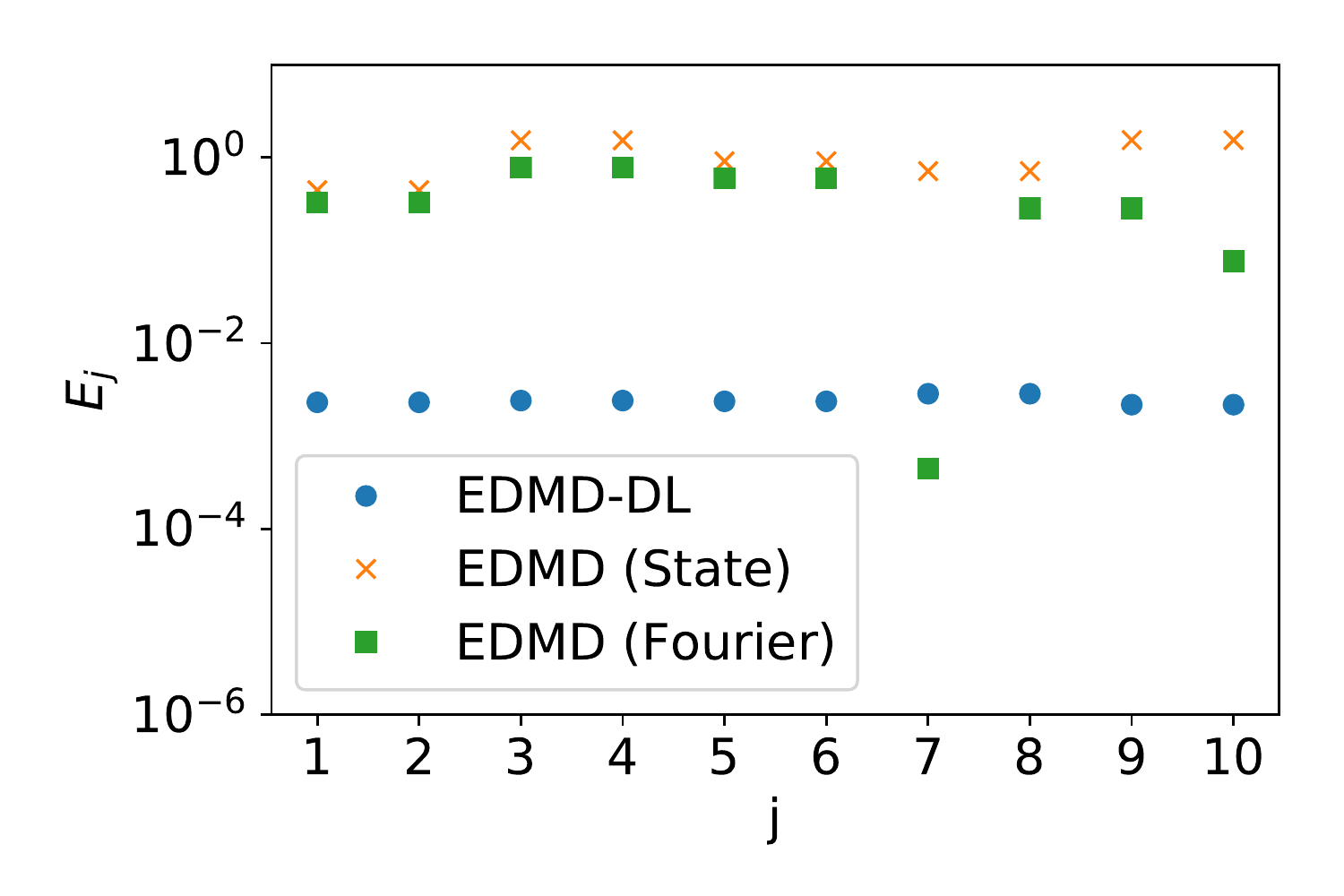}
	\end{center}
	\caption{Eigenfunction errors for various methods. Dictionary learning out-performs classical EDMD due to the data-adapted dictionary. Here, we used 25 dictionary elements for both EDMD-DL and EDMD with Hermite polynomial dictionaries, and 100 RBFs for EDMD with RBF dictionary.}
	\label{fig:kse_efuns}
\end{figure}
Next, we plot in Fig. \ref{fig:kse_traj} a reconstructed trajectory from a previously unseen initial condition. We see that classical EDMD with either choices of dictionaries cannot reproduce the fine-scale structures of the solution, whereas EDMD-DL achieves good reconstruction accuracy and manages to capture detailed behavior of the trajectory. Fig. \ref{fig:num_dict}(b) again confirms that EDMD-DL achieves good performance with smaller dictionary sizes. In fact, in this PDE case it is harder to pick a good dictionary and hence we see that the Fourier basis choice does not become better when more modes are included. This may be partially attributed to the fact that the Kuramoto-Sivashinsky PDE is known to possess an inertial manifold, so that the amplitudes of higher Fourier modes are effectively determined by those of the lower modes \cite{constantin2012integral,jolly1990approximate}. 

Lastly, in Fig. \ref{fig:kse_efuns} we observe that the eigenfunction errors $E_j$ (defined in \eqref{eq:Ej}) are much lower for EDMD-DL compared with classical EDMD. Here, instead of performing infinite-dimensional integration with some generic measure, we set $\rho$ to be the sample distribution of $u$ of the test trajectory used in Fig. \ref{fig:kse_traj}.

\section{Discussion\label{sec:discussion}}

Extended Dynamic Mode Decomposition approximates the spectrum of the Koopman operator, its eigenvalues, eigenfunctions, and modes. Until now, a dictionary in the form of a set of {\it a priori} chosen observables of the system states was not only necessary, but carefully choosing these was crucial to the performance of the method. In highly nonlinear and high-dimensional systems, such choices are hard to make. 
Our main contribution is formulating a problem \eqref{eq:dl_gen_min_prob} to find an optimal (in terms of the norm of the residual) choice of dictionary given the data. This allows for a low number of optimized dictionary functions to span a linear subspace on which the Koopman operator can be accurately approximated. To realize this algorithmically, we introduced an iterative algorithm in combination with a general approximator, in the form of a neural network. This leads to much more accurate reconstructions of the Koopman operator spectrum with fewer (adapted) observation functions. Furthermore, we see that these adaptive descriptions usually have greater reconstruction accuracy over longer trajectory lengths, even those exceeding the length of the training trajectories (see Fig. \ref{fig:duffing_traj}). 

These desirable properties of the EDMD-DL method enable a greater range of applications of spectral analysis of nonlinear dynamical systems in general. For instance, since fewer dictionary elements are needed by EDMD-DL (see Fig. \ref{fig:num_dict}), it can be readily applied to obtain accurate reconstruction for high dimensional ODE systems or PDE systems. Moreover, this linearization technique is also useful in enabling control theory of linear systems (which is a much studied subject) to be applied to nonlinear dynamics \cite{kwakernaak1972linear}. 

The use of neural networks as dictionary approximators is also interesting on its own. Besides being a universal approximator, a neural network can also be built with certain invariance properties if, so desired. For example, for applications involving spatially homogeneous PDEs, it is natural to seek eigenfunctions that are translation-invariant. Such symmetries can be built into the neural networks by considering convolution layers as their main building blocks instead of the fully-connected layers considered in this paper. These convolution neural networks (CNNs) have been extensively used in image processing and classification \cite{krizhevsky2012imagenet}, and are likely to be highly effective in dealing with PDE systems with spatially homogeneous and local interaction terms. Moreover, CNNs are also useful in picking up multi-scale features, and hence using CNNs as dictionary approximators is also expected to be useful in dealing with systems with dynamics that have multiple length scales \cite{lecun1995convolutional,lecun2010convolutional,krizhevsky2012imagenet,lecun2015deep}.

\section{Conclusion and outlook}\label{sec:conclusion}

In this paper, we combine modern machine learning approaches with the EDMD algorithm for estimating spectral decompositions of the Koopman operator. This allows us to address an important shortcoming of the EDMD algorithm, namely the choice of a dictionary. In the EDMD-DL framework, we regard the dictionary itself as an additional optimization variable. Consequently, we can seek the optimal finite-dimensional approximation of the Koopman operator given the size of the dictionary. This allows the application of the Koopman operator framework to a broader range of problems with improved reconstruction accuracies. 

There are many directions for future research. From the algorithmic point of view, the conditions which guarantee the convergence of Alg. \ref{alg:EDMDdict} can be studied. One can also explore variants of the algorithm with e.g. a different regularization term or a different function approximator that may be more suited for solving specific problems. 
If the data are known, or suspected, to live on a low-dimensional manifold, then the relation of the number of dictionary elements found with the number of ``generic observables" suggested by the Whitney, Nash and Takens embedding theorems \cite{sauer1991embedology} should be both interesting and informative to explore. 
The stochastic counter-part to the Koopman operator is the backward Kolmogorov operator for stochastic dynamics. It will be also interesting to apply this method to obtain spectral analysis of the evolution of expected values of observables driven by stochastic dynamics.

%For instance, an accurate spectral decomposition may allow one to compare the spectral properties of Koopman operators corresponding to pairs of dynamical systems that are spatial isomorphisms of each other. It was shown by von Neumann that spectral isomorphism of Koopman operators implies spatial isomorphism of dynamical systems under restricted conditions \cite{neumann1932operatorenmethode,halmos1942operator} (The general result is false. See \cite{halmos1957introduction}). With a robust numerical algorithm to compute Koopman spectra such as EDMD-DL, we can attempt to search for mappings between spectra of Koopman operators of two different dynamical systems, in order to gain insights into possible spatial mappings between dynamical systems. This may serve as a means to reduce complex dynamical systems to simpler ones through learned spatial mappings. 

\section*{Acknowledgements}
The work of IGK was partially supported by DARPA-MoDyL (HR0011-16-C-0116) and by the US National Science Foundation (ECCS-1462241).
IGK and FD are grateful for the hospitality and support of the IAS-TUM. FD is also grateful for the support from the TopMath Graduate Center of TUM Graduate School at the Technical University of Munich, Germany and from the TopMath Program at the Elite Network of Bavaria. EMB thanks the Army Research Office (N68164-EG) and the Office of Naval Research (N00014-15-1-2093). QL is grateful for the support of the Agency for Science, Technology and Research, Singapore.

\nocite{*}
\bibliography{ref}% Produces the bibliography via BibTeX.

\end{document}